\documentclass{amsart}
\usepackage[margin=1.5in]{geometry}

\usepackage[all]{xy}
\usepackage{tikz-cd}
\usepackage{enumerate}
\usepackage{amsfonts}
\usepackage{amssymb, amscd}
\usepackage{graphicx}
\usepackage{mathtools}
\usepackage{xcolor}
\usepackage{float}

\newcommand{\br}{\mathbf{r}}
\newcommand{\E}{\mathcal{E}}
\newcommand{\D}{\mathbb{D}}
\newcommand{\G}{\mathbb{G}}
\newcommand{\T}{\mathbb{T}}

\newcommand{\B}{\mathbb{B}}

\newcommand{\R}{\mathbb{R}}

\newtheorem{theorem}{Theorem}

\newtheorem{lemma}{Lemma}

\theoremstyle{definition}

\begin{document} 

\title{Reflections of convex bodies and their sections}  
\author{Jorge L. Arocha, Javier Bracho and Luis Montejano}

\begin{abstract}
The purpose of this paper is to study the reflections of a convex body. In particular, we are interested in orthogonal reflections of its sections that can be extended to reflections of the whole body. For this reason, we need  to study the case of a non-spherical ellipsoid, where a surprising structure arises (Section 2). These results allow us to give a new characterization of ellipsoids in terms of their reflections and, on the other hand, to prove a result deeply related to a conjecture due to K. Bezdek.
\end{abstract}
\maketitle

\section{Introduction}

Every nonempty planar section of a 3-dimensional ellipsoid is an ellipse. Hence, it admits at least two orthogonal reflections because the axes of the ellipse are reflection lines. 
Likewise, every  nonempty planar section of a 3-dimensional body of revolution admits a reflection line (which is coplanar with the axis of revolution).  Karoly Bezdek conjectured in \cite{Bezdek} that if every nonempty planar section of a $3$-dimensional convex body admits an orthogonal reflection then it must be either an ellipsoid or a body of revolution. 
Some work has been done in the direction of, or inspired by, this conjecture: 

Given a convex body $K$, let us define a \emph{Bezdek} plane  (for $K$, it is understood) to be a plane whose section with $K$ admits an orthogonal reflection; observe that non-intersecting planes are Bezdek planes.
We know that if $K$ is a convex body for which every plane is a Bezdek plane (the hypothesis of the conjecture) then, 
 through every interior point of $K$ there is a section which is a disk, \cite{M}. Moreover, if every plane is a Bezdek plane with a reflection line coplanar to a fixed line, then  
$K$ is a  body of revolution, \cite {JMM}.  For more related problems see \cite{R}.

Here we prove a weaker form of the conjecture. It  leads to the proof of a similar statement in higher dimensions.  First, we need to establish some terminology.

A  convex body $K\subset \R^n$ is an \emph{ellipsoid} if it is the image of a ball under an affine automorphism. We say that a convex body $K\subset \R^n$ is a $\mathbb S^1$-\emph{rotational body} if there is a $(n-2)$-dimensional flat $L$, called the \emph{axis of rotation} of $K$, with the property that for every plane $P$ orthogonal to $L$: $P\cap K$  is either empty, a single point  or a disk centered at $L\cap P$.  Clearly, in $\R^3$ an $\mathbb S^1$-rotational body is a body of revolution. 

Let $K$ be a convex body in $\R^n$ and let $\Gamma$ be a hyperplane.  Generalizing directly, we say that $\Gamma$ is a \emph{Bezdek} hyperplane if $\Gamma \cap K$ admits an orthogonal reflection. If, in addition, this orthogonal reflection in $\Gamma$ extends to a reflection
 $\bf{r}$ of $K$ we say that $\Gamma$ is \emph{strong-Bezdek} and we call $\bf{r}$ the \emph{strong reflection} associated to $\Gamma$. Observe that in $\R^3$, every plane is strong-Bezdek for an ellipsoid or a body of revolution. 

Our main purpose is to prove the following theorem which,  in dimension $3$, yields a weaker version of Bezdek's Conjecture. The other variation is on the set of planes considered.

\begin{theorem}\label{teo:main}
Let $K$ be a $n$-dimensional convex body.
If there exists an open subset of the  Grasmannian space whose elements are strong-Bezdek  hyperplanes, then $K$ is either an ellipsoid or a $\mathbb S^1$-rotational body. 
\end{theorem}

Given a reflection $\br : \R^n \to \R^n$, we must distinguish the following associated elements:  the \emph{mirror} is the hyperplane which is pointwise fixed, and the \emph{direction} which is a line  through the origin all of whose parallel lines are rigidly inverted at their intersection with the mirror. 
The direction line is complementary to the mirror hyperplane. And if the direction is orthogonal to the mirror, the reflection is said to be \emph{orthogonal}.  

We say that $\Gamma$ is a \emph{ground} hyperplane for the reflection $\br$ if $\Gamma$ contains the direction of $\br$ (which implies that $\Gamma$ is invariant under $\br$) and $\br$ restricted to $\Gamma$ 
is an orthogonal reflection in $\Gamma$.  
Observe that if $\Gamma$ is a ground hyperplane of a reflection, so is any parallel hyperplane, so that  a ground hyperplane may be assumed to pass through the origin.

Let $\br$ be a reflection  with direction $\ell$ and mirror $M$. If $\ell$ is not orthogonal to $M$ then the hyperplane $\Gamma$  generated by $ M \cap \ell^\perp$  and $\ell$ is the unique ground plane for $\br$.  On the other hand, if $\br$ is an orthogonal reflection then all the hyperplanes through $\ell$ are ground hyperplanes for $\br$.

Now, we can rephrase the definition of strong-Bezdek hyperplanes. First, a convex body $K\subset\R^n$ has to be explicit in the context. Then, a hyperplane $\Gamma$ is strong-Bezdek if it is a ground hyperplane of a \emph{reflection of} $K$ (i.e., that leaves $K$ invariant). If this is the case, all the parallel hyperplanes to $\Gamma$ are also strong-Bezdek because they are ground hyperplanes of the same strong reflection associated to $\Gamma$. So the notion of strong-Bezdek is defined up to parallelism. However, observe that even if all the parallel hyperplanes to $\Gamma$ are Bezdek, for $\Gamma$ to be strong-Bezdek one still needs that they have mirrors that make them Bezdek and align inside a hyperplane.

For the proof of Theorem \ref{teo:main}, we need to understand the question of how many reflections of a $n$-dimensional body $K$ are enough to prove that $K$ is an ellipsoid?  A possible answer  to this question is given in Theorem 2 of the next section.  

\section{Characterizing ellipsoids by means of reflections}

The main result of this section is the following characterization of ellipsoids, which is critical for our proof of Theorem~1. To state it simply, let us define a subset of a topological space to be \emph{thick} if it contains a nonempty open set; that is, if it has non-empty interior. Let $G_{1,n}$ and $G_{n-1,n}$ be the Grassmannian of lines and hyperplanes, respectively, through the origin in $\R^n$.

\begin{theorem}\label{teo:ClasifEllips}
A convex body $K\subset \mathbb R^n$ is an ellipsoid if one of the following holds:
\begin{enumerate}
\item[i)] the lines $\ell$ whose midpoints of all chords of $K$ parallel to $\ell$ lie in a hyperplane, are thick in $G_{1,n}$.
\item[ii)] the lines $\ell$ such that there is a reflection of $K$ with direction $\ell$, are thick in $G_{1,n}$.
\item[iii)] the ground hyperplanes of non orthogonal reflections of $K$, are thick in $G_{n-1,n}$.
\end{enumerate} 
\end{theorem}

First we need to study ellipsoids with special attention to the ground planes of their non orthogonal reflections. Then we use that knowledge to prove Theorem~\ref{teo:ClasifEllips}.

Let $\E\subset\R^n$ be an ellipsoid; we may assume it is centered at the origin. 
Recall that $G_{1,n}=G_1(\R^n)$ is the Grassmannian of lines through the origin in $\R^n$. Any line $\ell\in G_{1,n}$ is the direction of a unique reflection $\br_\ell$ of $\E$. Indeed, the midpoints of all non-trivial chords in $\E$ parallel to $\ell$ generate its \emph{mirror} $M(\ell)=M_\E(\ell)$, which is also the parallel at the origin of the tangent hyperplanes to $\E$ at the endpoints of  the \emph{chord} $\ell\cap\E$. So that $\br_\ell$ is an orthogonal reflection if and only if $\ell$ is a \emph{binormal} line to $\E$, that is, a line $\ell$ whose tangent hyperplanes at the endpoints of $\ell\cap\E$ are orthogonal to $\ell$; and this happens if and only if $M(\ell)$ is the orthogonal hyperplane to $\ell$.  These lines are also called \emph{axes} of the ellipsoid. 
 
 We will use that the bijective assignment of mirror hyperplanes to lines given by the ellipsoid $\E$ reverses inclusion; that is, given lines $\ell$ and $\ell^\prime$ in $G_{1,n}$, then
\begin{equation}\label{eq:polaridad}
\ell^\prime\subset M(\ell)\quad\iff\quad \ell\subset M(\ell^\prime)\,.
\end{equation} 
To see this, consider a quadratic form that defines $\E$. It has an associated symmetric bilinear form. What determines if a vector lies in the mirror of the line generated by another vector, is that the symmetric bilinear form vanishes on that pair. Then, the assertion (\ref{eq:polaridad}) easily follows.

Let $\B_\E\subset G_{1,n}$ be the set of binormal lines of $\E$, and denote its complement by 
\begin{equation}\label{diagonals}
\D_\E=G_{1,n}\setminus \B_\E\,;
\end{equation} 
 it is the set of directions of non orthogonal reflections of $\E$; we refer to the elements of $\D_\E$ as \emph{diagonal} lines of $\E$.  Assume henceforth that $\E$ is not a sphere so that $\D_\E$ is an open dense subset of  $G_{1,n}$, instead of being empty. 
 
 Denote by $G_{n-1,n}$, the Grasmannian space of hyperplanes through the origin in $\R^n$.
Our main goal is to understand the \emph{ground map} 
$$\gamma:\D_\E\to G_{n-1,n}$$
which sends $\ell\in\D_\E$ to the unique ground hyperplane $\gamma(\ell)=\Gamma_\ell\in G_{n-1,n}$ of the reflection $\br_\ell$ of $\E$ in the direction $\ell$. The binormals $\B_\E$ had to be removed from $G_{1,n}$ to have the ground map well defined as a function.
 
 Note that  $\gamma(\ell)$ varies continuously with $\ell$ because we can write
$$\gamma(\ell)=(M(\ell)\cap\ell^\perp)\oplus\ell\,,$$
where ``$\,^\perp$'' denotes the orthogonal complement and ``$\oplus$'' the sum of complementary subspaces of $\R^n$.
And note also that $\ell\subset\gamma(\ell)$ is an axis of the ellipsoid 
$\gamma(\ell)\cap \E$, because the reflection $\br_\ell$ restricted to the hyperplane $\gamma(\ell)$  is an  orthogonal reflection. 

Indeed, given a diagonal line $\ell\in \D_\E$, $\gamma(\ell)$ is the unique hyperplane $\Gamma$ containing $\ell$ for which $\ell$ is a binormal of $\Gamma \cap \E$.  To  see this, observe that for any $\ell\in G_{1,n}$ and $\Gamma\in G_{n-1,n}$ such that $\ell\subset\Gamma$, the mirror of $\ell$ with respect to the ellipsoid $\Gamma\cap\E$ is the intersection of $\Gamma$ with its mirror with respect to $\E$, which we can write
$$M_{\Gamma\cap\E}(\ell)=\Gamma\cap M(\ell)\,.$$
So that a line which is an axis of two distinct hyperplane sections of $\E$ has to be a binormal line of $\E$, and by definition it is not a diagonal. 

We can summarize these last remarks with the following:
 
\begin{lemma}\label{lem:imagen-inversa}
Let $\Gamma\in G_{n-1,n}$, then $\gamma^{-1}(\Gamma)$ consists of the axes of the ellipsoid $\Gamma\cap \E$ which are not axes of $\E$. \qed
\end{lemma}

\begin{theorem}\label{prop:CoveringMap}
There are open dense sets $\G\subset\D_\E$ and $\T\subset G_{n-1,n}$ 
such that the ground map of $\E$ restricts to a covering map $\gamma:\G\to \T$.
\end{theorem}

\proof  Given a line $\ell\subset\R^n$, let us denote by $\lambda(\ell)$ the \emph{length} of $\ell$, that is,
the length of the segment $\ell\cap\E$. Let 
$$\B_\E=\B_1 \cup \dots \cup \B_k$$
be the partition of the binormal lines of $\E$ given by their \emph{lengths}; that is, two binormals of $\E$,  belong to the same class $\B_i$ if and only if they have the same length $\lambda_i$. 
Denote by $V_i\subset \R^n$ the subspace generated by the lines in $\B_i\subset G_{1,n}$; so that $\B_i=G_1(V_i)$. We then have that 
$$\R^n=V_1\oplus\dots\oplus V_k$$ 
is the orthogonal partition of $\R^n$ by the eigenspaces of the symmetric matrix associated to a quadratic form that defines $\E$; so, we may refer to the $V_i$ as the \emph{eigenspaces} of $\E$. Observe that the partitions are non-trivial ($k\geq 2$) because we are assuming that $\E$ is not a sphere. 

Let us define the \emph{generic} diagonals
\begin{equation}\label{GenericLines}
\G=\{\,\ell\in G_{1,n} \,\mid\,\ell\not\subset V_i^\perp \,\, \text{and} \,\, \lambda(\ell)\neq\lambda_i\,, i=1,\dots,k\,\}\,.
\end{equation}
Observe that $\G\subset \D_\E$ because if $\ell\not\in\D_\E$ then $\ell\subset V_j$ for some $j$ and $V_j\subset V_i^\perp$ for all $i\neq j$. So that the ground map $\gamma$ is well defined on $\G$. 

Define the \emph{generic} hyperplanes to be
\begin{equation}\label{GenericHyperplanes}
\T=\{\,\Gamma\in G_{n-1,n}\,\mid\, V_i\not\subset\Gamma , i=1,\dots,k \,\}\,.
\end{equation}

First note that $\G\subset G_{1,n}$ and $\T\subset G_{n-1,n}$ are open dense subsets. We will prove that the ground map restricted to $\G$ is a covering map of $k-1$ sheets over $\T$, and in particular that
$$\gamma(\G)=\T.$$

To see that $\gamma(\G)\subset\T$
consider $\ell\in \G$, and suppose $\Gamma=\gamma(\ell)\notin \T$. This implies that $V_i\subset \Gamma$ for some $i$.
Since any two binormals of an ellipsoid having different lengths are orthogonal and both $\ell$ and the lines of $V_i$ are binormals of the ellipsoid $\Gamma\cap \E$, then $\lambda(\ell)\neq\lambda_i$ implies that $\ell \subset V_i^\perp\cap\Gamma\subset V_i^\perp$, which contradicts $\ell\in \G$.

Now consider a fixed hyperplane $\Gamma \in \T$. For each $i=1,\dots,k$, since $V_i\not\subset\Gamma$, we can choose a line $e_i\subset V_i$ which is not in $\Gamma$: if $V_i$ is a line, then $e_i=V_i$ and if $V_i$ is more than a line we choose $e_i\subset V_i$ to be the orthogonal line to $\Gamma\cap V_i$. Let us assume that $e_1,\dots,e_k$ correspond to the first coordinate axes of $\R^n$ so that we can refer to the subspace they generate as $\R^k$, but observe, for later on, that it really depends continuously on $\Gamma$. By construction, the binormal lines of $\E$ inside $\R^k$ are exactly the coordinate axis $e_1,\dots,e_k$ and they all have different lengths. We will now prove that in $\R^k$ there are exactly $k-1$ generic lines of different lengths which are the axes of $\Gamma\cap \E$. This implies, by Lemma~\ref{lem:imagen-inversa}, that 
$\# \gamma^{-1}(\Gamma)=k-1$ because generic lines are not binormal to $ \E $.

Let $\ell\subset \Gamma\cap\R^k$ be a binormal line to $\Gamma\cap\E$, that is, such that $M_{\Gamma\cap\E}(\ell)=M(\ell)\cap\Gamma$ is orthogonal to $\ell$. To prove that $\ell$ is generic, we first see that $\ell\not\subset V_i^\perp$. For this, it suffices to prove that $\ell\not\subset e_i^\perp$ because $V_i^\perp\subset e_i^\perp$. So suppose $\ell\subset e_i^\perp$. Since $e_i^\perp=M(e_i)$ ($e_i$ is binormal to $\E$) then, (\ref{eq:polaridad}) implies that $e_i\subset M(\ell)$. On the other hand, $e_i\not\subset\Gamma$ yields by simple dimension arguments that    
$$M(\ell)=(M(\ell)\cap\Gamma)\oplus e_i\,.$$
So that $M(\ell)$ is orthogonal to $\ell$ because both summands are orthogonal to $\ell$.
Then $\ell$ is binormal to $\E$, but being inside $e_i^\perp\cap\R^k$, $\ell$ has to be $e_j$ for some $j\neq i$, which contradicts that $e_j\not\subset\Gamma$. So, we have that $\ell\not\subset V_i^\perp$ for $i=1,\dots,k$.

We are left to prove $\lambda(\ell)\neq \lambda_i$ to conclude that $\ell\in\G$. Suppose not; that is, assume $\lambda(\ell)= \lambda_i$. Consider the plane $P=\ell\oplus e_i$ which intersects $\Gamma$ in $\ell$. The ellipse $P\cap\E$ has $e_i$ as an axis (because it is an axis of $\E$). Then if another chord through its center has the same length, it has to be a circle; because the length function of an ellipse reaches its unique extremes at the axes unless it is constant. This implies that the orthogonal line to $\ell$ in the plane $P$, $\ell'=P\cap\ell^\perp$, is contained in $M(\ell)$ but it is not in $\Gamma$,  so that $M(\ell)=(M(\ell)\cap\Gamma)\oplus\ell'$ is orthogonal to $\ell$. Then, $\ell$ is a binormal of $\E$, which implies, because of its length, that $\ell\subset V_i$; but this contradicts that $\ell\in\R^k$. 

So, we have proved that any axis of $(\Gamma\cap\E)\cap\R^k$ is a generic line.

Now consider two different lines $\ell$ and $\ell^\prime$ in $\Gamma\cap\R^k$ which are binormal to $\Gamma\cap\E$. Suppose they have the same length $\lambda$. Let $P=\ell\oplus\ell^\prime\subset\Gamma\cap\R^k$ be the plane generated by them. All of the lines in $G_1(P)$ are binormal to $\Gamma\cap\E$ and have length $\lambda$. This implies, in particular, that $\ell^{\prime\prime}=P\cap e_1^\perp\subset\R^k$ is a binormal line of $\Gamma\cap\E$. But  $\ell^{\prime\prime}$ is not generic (it is contained in $V_1^\perp$),  contradicting the conclusion of the previous paragraph. And therefore, $\lambda(\ell)\neq\lambda(\ell^\prime)$.

Summarizing, given $\Gamma\in\T$, since $\Gamma\cap\E\cap\R^k$ is an ellipsoid of dimension $k-1$, it has $k-1$ mutually orthogonal axes. We have proved that they are generic lines with different lengths, so that, by Lemma~\ref{lem:imagen-inversa}, $\gamma^{-1}(\Gamma)$ consists of exactly $k-1$ isolated points in $\G$ because all the other binormal lines to $\Gamma\cap\E$ are also binormal lines to $\E$.

Finally note that in a small neighborhood of $\Gamma$ the order of the lengths of the axes of the corresponding sections of $\E$ is well determined, which implies   not only that $\#\gamma^{-1}(\Gamma)= k-1$ but also that  $\gamma:\G\to \T$ is a $(k-1)$-covering space, as we wished. 
\qed 

Now, we can prove Theorem~\ref{teo:ClasifEllips}. It will ease reading to restate it:

\setcounter{theorem}{1}
\begin{theorem}
A convex body $K\subset \mathbb R^n$ is an ellipsoid if one of the following holds:
\begin{enumerate}
\item[i)] the lines $\ell$ whose midpoints of all chords of $K$ parallel to $\ell$ lie in a hyperplane, are thick in $G_{1,n}$.
\item[ii)] the lines $\ell$ such that there is a reflection of $K$ with direction $\ell$, are thick in $G_{1,n}$.
\item[iii)] the ground hyperplanes of non orthogonal reflections of $K$, are thick in $G_{n-1,n}$.
\end{enumerate} 
\end{theorem}

\proof 
With hypothesis (i) this theorem is a generalization of Brunn's characterization of ellipsoids (see Theorem 2.12.1 of  \cite{MMO} and \cite{ABM}), in the sense that only thickness of such lines is asked instead of being all of them. 

First observe that hypothesis (i) and (ii) are equivalent because there is a reflection of $K$ in the direction $\ell$ if and only if the midpoints of all the chords of $K$ parallel to $\ell$ lie in a hyperplane (the mirror). So, it is enough to consider (ii).

Let $\E$ be the unique ellipsoid of minimal volume containing $K$, see e.g. \cite{ABM}. Observe that any affine symmetry of $K$ is also an affine symmetry of $\E$. Indeed, if $g:\R^n\to\R^n$ is an affine map such that $K=g(K)$, then $g(\E)$ is an ellipsoid that contains $K$, but since $g$ preserves volume ($K$ is a body), the uniqueness of $\E$ implies that $\E=g(\E)$.

Let $\Lambda_K\subset G_{1,n}$ be the set of lines which are the direction of a reflection of $K$. As in the beginning of this section, let $\br_\ell$ be the reflection of $\E$ with direction $\ell$, so that $\br_\ell$ is well defined for all $\ell\in G_{1,n}$. But observe that only for $\ell\in \Lambda_K$ we have that $\br_\ell$ is also a reflection of $K$ (that is, $\ell\in\Lambda_k\iff\br_\ell(K)=K$).

Suppose (ii) is true, that is, assume $\Lambda_K$ is thick. We will prove $K=\E$.

To see that $K=\E$,  it is enough to prove that $K\cap\partial \E$ is open in $\partial \E$. Indeed, since $K\cap\partial \E=\partial K\cap\partial \E$ is non-empty and closed in $\partial \E$, if it is also open it has to be all of $\partial\E$, which implies that $K=\E$ because both sets are convex. 

Consider a point $p\in K\cap\partial\E$, and let $T_p$ be the tangent plane to $\E$ at $p$ (which is also a support plane to $K$ at $p$). Denote by $W_p\subset G_{1,n}$ the set of lines not parallel to the hyperplane $T_p$; it is an open dense set  of $G_{1,n}$ homeomorphic to $\R^{n-1}$, also known as an \emph{affine chart} of $G_{1,n}$.
For any line $\ell\in W_p$, we have that $\br_\ell(p) \neq p$ because the parallel to $\ell$ through $p$, which we will denote $\ell(p)$, intersects $\partial \E$ in another point which is  $\br_\ell(p)$; that is, $\ell(p)\cap\partial\E=\{p,\br_\ell(p)\}$. Since $\Lambda_K$ is thick in $G_{1,n}$, there clearly exists an open set $V\subset \Lambda_K\cap W_p$. Consider the set
$$V_p = \{\,\br_\ell(p)\,|\,\ell\in V\,\}\,.$$
It is contained in $K\cap\partial\E$ because $p\in K\cap\partial\E$ and $V\subset\Lambda_K$. And it is homeomorphic to $V$ because it can be identified with its image under the \emph{stereographic projection from $p$} that sends a line $\ell\in W_p$ to the point different from $p$ in $\ell(p)\cap\partial\E$, and which is a homeomorphism between $W_p$ and $\partial\E\setminus \{p\}$. So that $V_p$ is an open set in $\partial\E$. Finally, fix any line $\ell_0\in V\subset\Lambda_K$. Since $\br_{\ell_0}$ restricts to a homeomorphism of $K$ and of $\partial\E$, we have that $\br_{\ell_0}(V_p)\subset K\cap\partial\E$ is an open neighborhood of $p$ in $\partial\E$ and therefore $K\cap\partial\E$ is open in $\partial\E$: which proves $K=\E$. 

Next, assume (iii) and we will prove $K=\E$ showing that (ii) holds. Let $\Omega_K\subset G_{n-1,n}$ be the set of hyperplanes that are ground hyperplanes of non orthogonal reflections of $K$ (hypothesis (iii) is that $\Omega_K$ is thick). That it is non-empty, implies that $\E$ is not a sphere and so Theorem~\ref{prop:CoveringMap} with its notation applies. Since  $\Omega_K$ is thick, there exists a set $U \subset \Omega_K\cap\T$ which is open in $G_{n-1,n}$ and evenly covered by the ground map $\gamma$, so that
$$\gamma^{-1}(U)=U_1\cup\cdots\cup U_{k-1}$$
with each $U_i$ homeomorphic to $U$ via $\gamma$. Consider $\Gamma\in U$, the $k-1$ reflexions $\br_\ell$ with $\ell\in\gamma^{-1}(\Gamma)$ are all the non orthogonal reflexions of $\E$ with ground hyperplane $\Gamma$. At least one of them is a reflexion of $K$ because $\Gamma\in\Omega_K$ so that at least one of the lines in $\gamma^{-1}(\Gamma)$ is also in $\Lambda_K$. This implies that
$$U=\gamma(U_1\cap\Lambda_K)\cup\cdots\cup\gamma(U_{k-1}\cap\Lambda_K)\,.$$
Observe that $\Lambda_K$ is closed, so that by Baire's category theorem at least one of the $U_i\cap\Lambda_K$ has non empty interior.  This implies that $\Lambda_K$ is thick and therefore, by (ii), that $K$ is an ellipsoid. \qed

\section{Proof of the main theorem}

For the proof of Theorem 1, we need first the following lemma:
 \begin{lemma}\label{lem:NotRot} 
Let $K\subset \mathbb R^3$ be a convex body which is not $\mathbb S^1$-rotational. Then the set of orthogonal reflections of $K$ is finite
\end{lemma}

\proof  Let $G$ be the group of isometries of $K$. We may assume that $G\subset O_n$. Note that $K$ is an $\mathbb S^1$-rotational body if and only if the group $G$ contains a copy of the group $O_2$ of isometries of $\mathbb S^1\subset\R^2$. Therefore, $G$ does not contain a copy of the group $SO_2\subset O_2$ and, since it is a compact subgroup of $O_n$, it is a finite group. This implies that the set of orthogonal reflections of $K$ is finite.
\qed
\bigskip

\noindent \emph{Proof of Theorem 1}. 
Recall that a hyperplane $\Gamma\in G_{n-1,n}$ is strong-Bezdek if it is the ground plane of a reflection of $K$, and by assumption there is an open subset $U\subset G_{n-1,n}$ of such hyperplanes.

Suppose $K$ is not a $\mathbb S^1$-rotational body. By Lemma~\ref{lem:NotRot}, the set of orthogonal reflections of $K$ is finite. For each of them, the set of ground hyperplanes is homeomorphic to $\R\mathbb{P}^{n-2}\cong G_{n-2,n-1}\subset G_{n-1,n}$ (the hyperplanes that contain the direction of the reflection). So that the set of ground hyperplanes of orthogonal reflections of $K$ is, at most, a finite collection of flats $\R\mathbb{P}^{n-2}$ contained in  $\R\mathbb{P}^{n-1}\cong G_{n-1,n}$. In its complement, there clearly is an open set $U^\prime\subset U$ of ground hyperplanes of non orthogonal reflections of $K$. Therefore, by item (iii) of Theorem~\ref{teo:ClasifEllips}, $K$ is an ellipsoid. \qed

\bigskip


\end{document}